\documentclass{amsart}
\usepackage{amsmath}

\newtheorem{theorem}{\textbf{Theorem}}
\newtheorem{question}{\textbf{Question}}

\newtheorem{remark}{\textbf{Remark}}

\def\a {\alpha}
\def\b {\beta}

\def\t {\tau}
\def\p {\partial}

\def\Q {\mathbb{Q}}
\def\R {\mathbb{R}}

\def\QQ {\overline{\Q}}
\def\C {\mathbb{C}}

\def\e {\epsilon}
\def\d {\delta}

\theoremstyle{remark}

\numberwithin{equation}{section}

\begin{document}

\title[ON A STRONGER VERSION OF A QUESTION PROPOSED BY K. MAHLER]{ON A STRONGER VERSION OF A QUESTION PROPOSED BY K. MAHLER}

\author{DIEGO MARQUES}
\address{DEPARTAMENTO DE MATEM\'{A}TICA, UNIVERSIDADE DE BRAS\'ILIA, BRAS\'ILIA, DF, BRAZIL}
\email{diego@mat.unb.br}


\author{CARLOS GUSTAVO MOREIRA}
\address{INSTITUTO DE MATEM\' ATICA PURA E APLICADA, RIO DE JANEIRO, RJ, BRAZIL}
\email{gugu@impa.br}

%    General info


\begin{abstract}
In 1902, P. St\"{a}ckel proved the existence of a transcendental function $f(z)$, analytic in a neighbourhood of the origin, and with the property that both $f(z)$ and its inverse function assume, in this neighbourhood, algebraic values at all algebraic points. Based on this result, in 1976, K. Mahler raised the question of the existence of such functions which are analytic in $\mathbb{C}$. Recently, the authors answered positively this question. In this paper, we prove a much stronger version of this result by considering other subsets of $\C$.
\end{abstract}

\maketitle

%%  SECTION 1
\section{Introduction}

A \textit{transcendental function} is a function $f(x)$ such that the only complex polynomial satisfying $P(x, f(x)) =0$ for all $x$ in its domain, is the null polynomial. For instance, the trigonometric functions, the exponential function, and their inverses. 

The study of the arithmetic behavior of transcendental functions at complex points has attracted the attention of many mathematicians for decades. The first result concerning this subject goes back to 1884, when Lindemann proved that the transcendental function $e^z$ assumes transcendental values at all nonzero algebraic point. In 1886, Strauss tried to prove that an analytic transcendental function  cannot assume rational values at all rational points in its domain. However, in 1886, Weierstrass supplied him with a counter-example and also stated that there are transcendental entire functions which assume algebraic values at all algebraic points. This assertion was proved in 1895 by St\"{a}ckel \cite{19} who established a much more general result: for each countable subset $X\subseteq \C$ and each dense subset $Y\subseteq \C$, there exists a transcendental entire
function $f$ such that $f(X) \subseteq Y$. In another construction, St\"{a}ckel \cite{20} produced a transcendental function $f(z)$, analytic in a neighbourhood of the origin, and with the property that both $f(z)$ and its inverse function assume, in this neighbourhood, algebraic values at all algebraic points. Based on this result, in 1976, Mahler \cite[p. 53]{bookmahler} suggested the following question

\begin{question}\label{q1}
Does there exist a transcendental entire function
\[
f(z)=\sum_{n=0}^{\infty}a_nz^n,
\]
with rational coefficients $a_n$ and such that the image and the preimage of $\QQ$ under $f$ are subsets of $\QQ$?
\end{question}

We refer the reader to \cite{bookmahler,wal1} (and references therein) for more about this subject. In a recent work, the authors \cite{DG} answered positively this question by proving the existence of uncountable many of these functions. 

In this paper, we prove a result which generalizes the main theorem of \cite{DG}. More precisely, we have
\begin{theorem}\label{1}
Let $X$ and $Y$ be countable, dense and closed for complex conjugation subsets of $\C$. Suppose that either both $X\cap \R$ and $Y\cap \R$ are dense in $\R$ or both intersections are the empty set and that if $0\in X$, then $Y\cap \Q\neq\emptyset$. Then, there are uncountably many transcendental entire functions
\[
f(z)=\sum_{n=0}^{\infty}a_nz^n,
\]
with rational coefficients $a_n$ and such that $f(X)=Y$, $f^{-1}(Y)=X$ and $f'(\alpha)\ne 0$, for all $\alpha\in X$.
\end{theorem}

Let us describe in a few words the difficulty of this problem when compared with the authors' proof for the original Mahler's question. First, in \cite{DG} it was strongly used that the desired function $f$ is a limit of polynomials $f_n$  with real algebraic coefficients. In particular, in each step, $f_n(A)$ and $f_n^{-1}(A)$ are finite subsets of $\QQ$, for every finite subset $A\subseteq \QQ$. This holds, clearly, because $\QQ$ is an algebraically closed field. In the proof of the theorem in this work, the sets $X$ and $Y$ do not need to have this property (algebraic closeness) and so in each step we need to ``force"\ these desired requirements. For that, we use the implicit function theorem and the inverse function theorem.

\begin{remark}
We point out that Mahler, in his book, said that the Question \ref{q1} was ``unknown whether there exists an entire function of this kind where the coefficients may be arbitrary complex numbers". In fact, in this case, we can make the hypothesis of Theorem \ref{1} weaker. More precisely, after an easy adaptation of the proof of this theorem, we can prove that: \textit{Let $X$ and $Y$ be countable and dense subsets of $\C$. Then, there are uncountably many transcendental entire functions $f$ such that $f(X)=Y$, $f^{-1}(Y)=X$ and $f'(\alpha)\ne 0, \forall \alpha\in X$. Moreover, assuming that $Y\cap {\mathbb Q}[i]\ne 0$ in the case when $0\in X$, we may request that all coefficients in the Taylor representations of such functions $f$ belong to ${\mathbb Q}[i]$.}
\end{remark}

{\bf Acknowledgements:} We would like to thank Prof. Michel Waldschmidt for pointing out an important correction in a previous version of this work. The second author also would like to thank Prof. Nicolau Saldanha, who invented in collaboration with him a problem for the Brazilian Mathematical Olympiad for University Students of 2008 asking to prove that there is an entire function which extends an increasing diffeomorphism $f:\mathbb R \to \mathbb R$ such that $f(\mathbb Q)=\QQ$, whose ideas were useful in this work.

\section{The proof}

In order to simplify our presentation, we use the familiar notation $[a, b] = \{a, a + 1,\ldots, b\}$, for integers $a < b$.

In order to avoid unnecessary repetitions, we shall consider here only the case in which $X\cap \R$ and $Y\cap \R$ are dense in $\R$. Suppose, without loss of generality, that $0\in X$ and that $r\in Y\cap \Q$. Let $\{\a_1,\a_2,\a_3,\ldots\}$ be an enumeration of $X$ such that for any $n\geq 1$, the numbers $\a_{3n-1},\a_{3n}\notin \R$ with $\a_{3n}=\overline{\a}_{3n-1}$ and $\a_{3n+1}\in \R$. Also, let $\{\b_1,\b_2,\b_3,\ldots\}$ be an enumeration of $Y$ such that for any $n\geq 1$, the numbers $\b_{3n-1},\b_{3n}\notin \R$ with $\b_{3n}=\overline{\b}_{3n-1}$ and $\b_{3n+1}\in \R$ (take $\a_1=0$ and $\b_1=r$). Now, let us construct our desired function inductively. 

Define $f_1(z)=z+r$, then $f_1(\a_1)=r\in Y$ and $f_1^{-1}(\b_1)=0\in X$. Now, we want to construct a sequence of polynomials $f_2(z), f_3(z),\ldots$ recursively of the form
\[
f_{m+1}(z)=f_m(z)+z^mh_m(z)P_m(z),
\]
with
\begin{itemize}
\item[(i)] $P_m\in \R[z]$ and $f_m(z)=\sum_{i=0}^{t_m}a_iz^i$ with $t_m\ge m$;
\item[(ii)] Let $X_m=\{\a_2,\a_3\ldots, \a_{3m-2}\}$ and $\tilde X_m=\{\tau\in f^{-1}_m(\{\b_1,\ldots, \b_{3m-2}\})|\tau\ne 0, \tau\in X, f_m'(\tau)\neq 0\}$. Then $P_m(z)=\prod_{\tau \in X_m \cup \tilde X_m}(x-\tau)^2$. Moreover, $P_m(z)\mid P_{m+1}(z), \forall m\ge 1$;
\item[(iii)]$f_m(\tau)\in Y, f_m'(\tau)\ne 0, \forall \tau\in X_m$ and, for each integer $j\in [1, 3n-2]$, there is $\tau\in \tilde X_m$ such that $f_m(\tau)=\b_j$. 
\item[(iv)] $0<L(h_m P_m)<\nu_m:=\frac{1}{m^{m+2+\deg (h_m P_m)}}$. 
%\item[(iv)] $0<\max\{|\e_m|, |\d_m|\}<\nu_m:=\frac{1}{L(P_m)s_mm^{m+3+\deg P_m}}$. Here $s_m$ depends only on $m$ and $f_{m-1}^{-1}(\{\b_1,\ldots, \b_{3m-2}\})$;
\item[(v)] $a_0,a_1,\ldots, a_m\in \Q$ and $a_k\ne 0,$ for $k\in [1,m]$.
%\item[(vi)] $f_m(\{\a_1,\ldots, \a_{3m-2}\})\subseteq Y$;
\item[(vi)] $f^{-1}_{m+1}(\{\b_1,\ldots, \b_{3m-2}\})\cap B(0,r_m)=f^{-1}_m(\{\b_1,\ldots, \b_{3m-2}\})\cap B(0,r_m)=f^{-1}_m(\{\b_1,\ldots, \b_{3m-2}\})\cap \overline{B(0,r_m)}=\tilde X_m\cap B(0,r_m)\subseteq X$, for some suitable choice of $r_m$ with $m<r_m<m+1$;
%%\item[(viii)] $f_m'(\theta)\neq 0$, for all $\theta\in \{\a_1,\ldots, \a_{3m-2}\}\cup f^{-1}_m(\{\b_1,\ldots, \b_{3m-2}\})$.
\end{itemize}
Here $L(P)$ denotes the \textit{length} of the polynomial $P$ (the sum of the absolute values of its coefficients).

The polynomials $h_n$ have the form $\sum_{j=0}^{s_n}(\e_{n,j}+\d_{n,j}z)z^{1-\hat{\d}_{j,0}}\tilde P_{n,j}(z)$ (here $s_n\le \tilde s_n:=3+(3n+1)\max\{\deg f_n,n+1+\deg P_n\}$ is a natural number which will be chosen later), where $\hat{\d}_{j,0}$ is $1$ if $j=0$ and $0$ otherwise and $\tilde P_{n,j}$ are monic polynomials with real coefficients. 
%Here $s_m$ depends only on $m$ and $f_{m-1}^{-1}(\{\b_1,\ldots, \b_{3m-2}\})$;
The requested function will have the form $f(z)=\lim_{m\to \infty}f_m(x)$, where, for each $m\ge 1$, $f_{m+1}(z)=r+z+\sum_{1\le n\le m}\sum_{j=0}^{s_n}(\e_{n,j}+\d_{n,j}z)z^{n+2-\hat{\d}_{j,0}}P_{n,j}(z)$, where $P_{n,j}(z)=P_n(z)\tilde P_{n,j}(z)$. Thus we have 
$$f(z)=r+z+\sum_{n\geq 2}\sum_{j=0}^{s_n}(\e_{n,j}+\d_{n,j}z)z^{n+2-\hat{\d}_{j,0}}P_{n,j}(z).$$ 
We will have, for $0\le j<s_n$, $P_n(x)|P_{n,j}(x)|P_{n,j+1}(x)|P_{n+1}(x)$.

In each step, we shall choose $0<\max\{|\e_{n,j}|, |\d_{n,j}|\}<\nu_{n,j}:=\frac{1}{L(P_{n,j})\tilde s_nn^{n+3+\deg P_{n,j}}}$. Since $|P_{n,j}(z)|\leq L(P_{n,j})\max\{1,|z|\}^{\deg P_{n,j}}$, we have that for all $z$ belonging to the open ball $B(0,R)$
\[
\left|\sum_{j=0}^{s_n}(\e_{n,j}+\d_{n,j}z)P_{n,j}(z)\right|<\frac{R+1}{n}\left(\frac{\max\{1,R\}}{n}\right)^{n+2+\deg P_{n,j}}.
\]
Thus $f$ is an entire function, since the series $\sum_{n\geq 2}\sum_{j=0}^{s_n}(\e_{n,j}+\d_{n,j}z)z^{n+2-\hat{\d}_{j,0}}P_{n,j}(z)$ converges uniformly in any of these balls.

Suppose that we have a function $f_{n}$ satisfying (i)-(viii). Now, let us construct $f_{n+1}$ with the desired properties. 

Let $B_n=B(0,r_n)$. We will define $f_{n,0}(z)$ as
\[
f_{n,0}(z)=f_n(z)+\e_{n,0}z^{n+1}P_n(z),
\]
for some choice of a (small) real number $\e_{n,0}$ (so, in this case, we take $P_{n,0}(z)=P_n(z)$). 
We are assuming, by induction hypothesis, that $f^{-1}_n(\{\b_1,\ldots, \b_{3n-2}\})$ does not intersect $\partial B_n$, and that, for every $\tau\in f^{-1}_n(\{\b_1,\ldots, \b_{3n-2}\})\cap B_n$, we have $\tau\in X$ and $f_n'(\tau)\ne 0$ - it follows that $\tau$ is a double zero of $P_n$.  We will choose $\e_{n,0}\in \R$ such that
\[
0<\e_{n,0}<\min_{i\in [1,3n-2]}\frac{\min_{|z|=r_n}|f_n(z)-\b_i|}{\max_{|z|=r_n}|z^{n+1}P_n(z)|}.
\]
It follows, by Rouch\' e's theorem, that the number of zeros (counted with multiplicity) of $f_n(z)-\b_i$ and $f_{n,0}(z)-\b_i$ belonging to $B_n$ are equal. Since every zero of $f_n(z)-\b_i$ in $B_n$ is a zero of $f_{n,0}(z)-\b_i$ and every zero of $f_n(z)-\b_i$ in $B_n$ is simple, we have that $f_{n,0}^{-1}(\b_i)\cap B_n= f_{n}^{-1}(\b_i)\cap B_n$, for all $i\in [1,3n-2]$.

We will show that, except for a finite set of values of $\e_{n,0}$, for any $w\in \C$ with $f_{n,0}(w)\in \{\b_1,\ldots, \b_{3n+1}\}$ we have $f_{n,0}'(w)\ne 0$. If $w$ is a root of $P_n(z)$, then $f'_{n,0}(w)=f'_n(w)\neq 0$. Then, $w$ is a simple root of $f_n(z)-\b_i$. Otherwise, $P_n(w)\ne 0$. We have $f_{n,0}(z)=f_n(z)+\e_{n,0}g(z)$, where $g(z)=z^{n+1}P_n(z)$. If $f_{n,0}'(w)=0$, we should have $f_n(w)+\e_{n,0}g(w)=\b_i$ and $f_n'(w)+\e_{n,0}g'(w)=0$. Defining $h_i(z)=f_n(z)-\b_i$, we have $h_i'(z)=f_n'(z)$, and thus $h_i(w)+\e_{n,0}g(w)=0$ and $h_i'(w)+\e_{n,0}g'(w)=0$. Let now $\psi_i(z)=-h_i(z)/g(z)$. Since $h_i(w)+\e_{n,0}g(w)=0$, we have $\psi_i(w)=-h(w)/g(w)=\e_{n,0}$. Moreover, $\psi_i'(w)=\frac{h_i(w)g'(w)-h_i'(w)g(w)}{g(w)^2}=0$ since $h_i(w)g'(w)-h_i'(w)g(w)=(h_i(w)+\e_{n,0}g(w))g'(w)-(h_i'(w)+\e_{n,0}g'(w))g(w)=0$. This implies that $\e_{n,0}$ is a singular value of $\psi_i(w)$, and the set of singular values of a rational function is finite, which concludes the argument.

In fact, there is an interval $I_{n,0}=(0,\e_0)$ of possible values for $\e_{n,0}$ for which the above properties hold. Since $P_{n,0}(0)\neq 0$, then we can choose $\e_{n,0}\in I_{n,0}$ such that the coefficient of $z^{n+1}$ in $f_{n,0}(z)$ is a nonzero rational (since this coefficient is equal to the coefficient of $z^{n+1}$ in $f_n(z)$, which is a real number, added by $\e_{n,0}P_{n}(0)\in \R$). Since in all further perturbations we will add multiples of $z^{n+2}$, $a_{n+1}$ will be equal to this coefficient, and thus is a nonzero rational number.

Let $r_{n+1}$ be a real number with $n+1<r_{n+1}<n+2$ and such that the intersection
$f_{n,0}^{-1}(\{\b_1,\ldots, \b_{3n+1}\})\cap \partial B(0,r_{n+1})$ is the empty set. For simplicity, we write $B_{n+1}=B(0,r_{n+1})$. Let $f_{n,0}^{-1}(\{\b_1,\ldots, \b_{3n+1}\})=\{\t_1,\t_2,\dots,\t_{m_n}\}$. Notice that each element of $\{\b_1,\ldots, \b_{3n+1}\}$ has at least one (indeed $\deg f_{n,0}$) pre-images in $\{\t_1,\t_2,\dots,\t_{m_n}\}$. For $j\in [1, m_n]$, fix a small ball $B(\t_j,\eta_j)$ which does not intersect $\partial B_{n+1}$ in such a way that the balls $B(\t_j,\eta_j)$ are disjoint. The number $s_n$ of further steps in the construction of $f_{n+1}$ will be the number of complex conjugation classes in $(\{\a_{3n-1},\a_{3n},\a_{3n+1}\}\cup \{\t_1,\t_2,\dots,\t_{m_n}\})\setminus(X_n\cup \tilde X_n)$. The following $s_n$ perturbations $f_{n,j}$ (with $j\in [1, s_n]$) of $f_{n,0}$ will be taken so close to $f_{n,0}$ that, by Rouch\'e's theorem as before, for each $j\le s_n$, the number of zeros (counted with multiplicity) of $f_{n,0}(z)-\b_i$ and $f_{n,j}(z)-\b_i$ in the balls $B_{n+1}$, $B_n$ and $B(\t_j,\eta_j), j\le m_n$ are equal, for all $i\in [1,3n+1]$ (in order to guarantee that we will assume that the absolute values of the coefficients $\e_{n,j},\d_{n,j}$ are very small - we will then say that $\e_{n,i}, \d_{n,i}$ are admissible). Notice that in the balls $B(\t_j,\eta_j), j\le m_n$ these numbers are equal to one, since $\t_j$ is a simple zero of $f_{n,0}(z)-\b_i$ for some $i\in [1,3n+1]$. We will take $f_{n+1}=f_{n,s_n}$, and so the number of zeros (counted with multiplicity) of $f_{n,0}(z)-\b_i$ and $f_{n+1}(z)-\b_i$ in the balls $B_{n+1}$, $B_n$ and $B(\t_j,\eta_j), j\le m_n$ will be equal, for all $i\in [1,3n+1]$. In particular, as before, $f_{n+1}^{-1}(\b_i)\cap B_n= f_{n}^{-1}(\b_i)\cap B_n$, for all $i\in [1,3n-2]$ and, for $j\in [1, m_n]$, $f_{n+1}$ has only one zero in $B(\t_j,\eta_j)$, which is simple.

Now, we define $f_{n,1}(z)$ by  
\[
f_{n,1}(z)=f_{n,0}(z)+\e_{n,1}z^{n+2}P_n(z),
\]
(so $P_{n,1}(z)$ is also equal to $P_n(z)$). Note that $f'_{n,1}(y)\neq 0$, for all $y\in \{\a_1,\ldots,\a_{3n-2}\}\cup f_n^{-1}(\{\b_1,\ldots, \b_{3n+1}\})$. Let $I_{n,1}$ be the interval of the admissible $\e_{n,1}$'s. If $\a_{3n+1}\in \R$ is a root of $P_n(z)$, then $f_{n,1}(\a_{3n+1})=f_{n}(\a_{3n+1})\in Y$ (by the definition of the roots of $P_n(z)$ in (ii)). Hence, in this case, we are done. So, suppose that $\a_{3n+1}$ is not a zero of $P_{n,1}(z)$. Then, by density of $Y\cap \R$, there exists an $\e_{n,1}\in I_{n,1}$ such that
\[
f_{n,1}(\a_{3n+1})=f_{n,0}(\a_{3n+1})+\e_{n,1}\a_{3n+1}^{n+2}P_{n,1}(\a_{3n+1})\in Y\cap \R.
\]
Define $f_{n,2}(z)=f_{n,1}(z)+(\e_{n,2}+\d_{n,2}z)z^{n+2}P_{n,2}(z)$, where $P_{n,2}(z)=P_{n,1}(z)(z-\a_{3n+1})\in \R[x]$. By applying a similar argument as before, there is an interval of possible choices of $\e_{n,2}, \d_{n,2}$ which are admissible. If $P_{n,2}(\a_{3n-1})=0$, then as before we obtain $f_{n,2}(\a_{3n-1})\in Y$. So, let us suppose that $P_{n,2}(\a_{3n-1})\neq 0$. Since $\a_{3n-1}\not\in \R$, then $1,\a_{3n-1}$ are $\R$-linearly independent and thus for $\e_{n,2},\d_{n,2}$ varying in $I_{n,2}$, then $\e_{n,2}+\d_{n,2}\a_{3n-1}$ covers a paralelogram in $\R^2$ (in particular, with non-empty interior). Thus, since $Y$ is dense in $\C$, for a suitable choice of $\e_{n,2}$ and $\d_{n,2}$ belonging to $I_{n,2}$, we have that  
\[
f_{n,2}(\a_{3n-1})=f_{n,1}(\a_{3n-1})+(\e_{n,2}+\d_{n,2}\a_{3n-1})\a_{3n-1}^{n+2}P_{n,2}(\a_{3n-1})\in Y.
\]
Moreover, $f_{n,2}(\a_{3n})=f_{n,2}(\overline{\a}_{3n-1})=\overline{f_{n,2}(\a_{3n-1})}\in \overline{Y}=Y$.

Let us define $f_{n,3}(z)$ by
\[
f_{n,3}(z)=f_{n,2}(z)+\e_{n,3}z^{n+2}P_{n,3}(z),
\]
where $P_{n,3}(z)=P_{n,2}(z)(z-\a_{3n-1})(z-\a_{3n})\in \R[z]$. Set $g_{n,3}(z)=z^{n+2}P_{n,3}(z)$. Since we are supposing that $\a_{3n-i}$ is not a root of $P_{n,1}(z)$, for $i\in [-1,1]$, then $g'_{n,3}(\a_{3n-i})\neq 0$ and so unless of three exceptions, we can choose $\e_{n,3}$ admissible in $I_{n,3}$ such that $f'_{n,3}(\a_{3n-i})\neq 0$, for $i\in [-1,1]$.  

Let $J$ be the set of indices $j\in [1, m_n]$ such that $Im(\t_j)\ge 0$ and $\t_j$ does not belong to $\{\a_{3n-1},\a_{3n},\a_{3n+1}\}\cup X_n\cup \tilde X_n\cup \{0\}$. Let $J=\{j_1,j_2,\dots,j_{s_n-3}\}$. For each $i\in [1, s_n-3]$, we will do a perturbation as below. Let $z_i$ be the only element of $f_{n,2+i}^{-1}(\{\b_1,\ldots, \b_{3n+1}\})$ in $B(\t_{j_i},\eta_{j_i})$. In the $i$-th step we will guarantee that the element of $f_{n,3+i}^{-1}(\{\b_1,\ldots, \b_{3n+1}\})$ in $B(\t_{j_i},\eta_{j_i})$ will belong to $X$, and in the further steps the image of these element by the maps $f_{n,2+j}$ will remain unchanged.

Suppose that $z_1\not\in \R$ (the real case can be done in a similar and easier way, as we will see in the next case). Define $f_{n,4}(z)$ by
\[
f_{n,4}(z)=f_{n,3}(z)+(\e_{n,4}+\d_{n,4}z)z^{n+2}P_{n,4}(z),
\]
where $P_{n,4}(z)=P_{n,1}(z)(z-\a_{3n-1})^2(z-\a_{3n})^2(z-\a_{3n+1})^2$. To simplify, we write $\e:=\e_{n,4}$, $\d:=\d_{n,4}$ and $F(\e,\d,z)=f_{n,4}(z)$. Since $\frac{\p F}{\p z}(0,0,z_1)=f_{n,3}'(z_1)\neq 0$, then the implicit function theorem ensures the existence of an implicit function $G(\e,\d,z)$, defined in a neighbourhood of $(0,0,z_1)$, such that $F(\e,\d,G(\e,\d,w))=w$, where $z_1=G(0,0,w)$ (with $w\in \{\b_1,\ldots,\b_{3n+1}\}$). We want to show the existence of $\e$ and $\d$ admissible such that $G(\e,\d,w)\in X$. By using the chain rule to derive $G$ at $(0,0,w)$ implicitly in the relation $F(\e,\d,G(\e,\d,w))=w$, we obtain
\[
z_1^{n+2}P_{n,4}(z_1)+f'_{n,3}(z_1)\frac{\p G}{\p \e}(0,0,w)=0
\]
and
\[
z_1^{n+3}P_{n,4}(z_1)+f'_{n,3}(z_1)\frac{\p G}{\p \d}(0,0,w)=0.
\]
Thus,
\begin{center}
$\frac{\p G}{\p \e}(0,0,w)=-\frac{z_1^{n+2}P_{n,4}(z_1)}{f'_{n,3}(z_1)}$ and $\frac{\p G}{\p \d}(0,0,w)=-\frac{z_1^{n+3}P_{n,4}(z_1)}{f'_{n,3}(z_1)}$.
\end{center}
Here we used that $f_{n,3}'(z_1)\neq 0$. Since $P_{n,4}(z_1)\neq 0$ and $z_1\not\in \R$, then the numbers $\frac{\p G}{\p \e}(0,0,w)$ and $\frac{\p G}{\p \d}(0,0,w)$ are $\R$-linearly independent. Thus, the determinant of Jacobian of $G$ at $(0,0,w)$ is nonzero. Therefore, by the inverse function theorem, the function $G(\e,\d,w)$ is a diffeomorphism of a small neighbourhood of $(\e,\d)=(0,0)$ into a neighbourhood of $z_1$. Since $X$ is a dense set, then there exist $\e,\d$ admissible such that $G(\e,\d,w)\in X$. 

Now, we suppose that $z_2\in \mathbb R$ and define $f_{n,5}(z)$ by
\[
f_{n,5}(z)=f_{n,4}(z)+\e_{n,5}z^{n+2}P_{n,5}(z),
\]
where $P_{n,5}(z)=P_{n,4}(z)(z-z_1)^2(z-\overline{z_1})^2\in \R[z]$ (note that by construction $G(\e,\d,\overline{w})$ also belongs to $X$). A similar (and indeed simpler, since we have only one parameter $\e=\e_{n,5}$ instead of the two parameters $\e$, $\d$) argument as in the previous case implies the existence of $\e=\e_{n,5}$ admissible such that the only pre-image by $f_{n,5}$ of $f_{n,4}(z_2)=\b_i\in\mathbb R$ in $B(\t_i,\eta_i)$ belongs to $X$.

So, following this construction, we will obtain at the end a function $f_{n,s_n}$ such that $f_{n,s_n}^{-1}(\{\b_1,\ldots, \b_{3n+1}\})\cap (B_{n+1} \cup (\cup_{j=1}^{m_n}B(\t_j,\eta_j)))$ is a subset of $X$. We set $f_{n+1}(z):=f_{n,s_n}(z)$. This function will satisfy the items (i)-(vi) (with $n$ replaced by $n+1$). 

Thus, by construction, the function 
\[
f(z)=\lim_{m\to \infty}f_m(z)=r+z+\sum_{n\geq 2}\sum_{j=0}^{s_n}(\e_{n,j}+\d_{n,j}z)z^{n+2-\hat{\d}_{j,0}}P_{n,j}(z)=\sum_{n\geq 0} a_nz^n
\]
is entire, $f(X)=Y, f^{-1}(Y)=X$ and $a_n\in \Q$ as desired. Indeed, for each $m\ge 1$, $f_{k+1}(\a_m)=f_k(\a_m)\in Y$ for all $k$ such that $3k-2\ge m$ (in particular for all $k\ge m$). In particular since $f=\lim_{k\to\infty}f_k$, we have $f(\a_m)=f_m(\a_m)\in Y$. So $f(X)\subset Y$. On the other hand, since $f^{-1}_{m+1}(\{\b_1,\ldots, \b_{3m-2}\})\cap B(0,r_m)=f^{-1}_m(\{\b_1,\ldots, \b_{3m-2}\})\cap B(0,r_m)=f^{-1}_m(\{\b_1,\ldots, \b_{3m-2}\})\cap \overline{B(0,r_m)}=\tilde X_m\cap B(0,r_m)\subseteq X$ and $B(0,r_m)\subset B(0,r_{m+1}), \forall m\ge 1$, it follows by induction that, for each $k\ge m$, $f^{-1}_k(\{\b_1,\ldots, \b_{3m-2}\})\cap B(0,r_m)=f^{-1}_m(\{\b_1,\ldots, \b_{3m-2}\})\cap B(0,r_m)$. This implies that 
$$f^{-1}(\{\b_1,\ldots, \b_{3m-2}\})\cap B(0,r_m)=f^{-1}_m(\{\b_1,\ldots, \b_{3m-2}\})\cap B(0,r_m)\subseteq X.$$ 
Indeed, $f=\lim_{k\to\infty}f_k$, and so 
$$f^{-1}(\{\b_1,\ldots, \b_{3m-2}\})\cap B(0,r_m)\supset f^{-1}_m(\{\b_1,\ldots, \b_{3m-2}\})\cap B(0,r_m).$$ 
On the other hand, if there were another element $w$ of  $f^{-1}(\{\b_1,\ldots, \b_{3m-2}\})\cap B(0,r_m)$, it should be at positive distance of the finite set $f^{-1}_m(\{\b_1,\ldots, \b_{3m-2}\})\cap B(0,r_m)$, but, since $f=\lim_{n\to \infty} f_k$, arbitrarily close to $w$ there should be, for $k$ large, an element of $f^{-1}_k(\{\b_1,\ldots, \b_{3m-2}\})$ (again by Rouch\' e's theorem), which contradicts the equality $f^{-1}_k(\{\b_1,\ldots, \b_{3m-2}\})\cap B(0,r_m)=f^{-1}_m(\{\b_1,\ldots, \b_{3m-2}\})\cap B(0,r_m)$. So $f^{-1}(Y)\subset X$ (and so, since $f(X)\subset Y$, $f^{-1}(Y)=X$).
Moreover, since for each integer $j$ with $1\le j\le 3n-2$, there is $\tau\in \tilde X_n$ such that $f_n(\tau)=\b_j$, and, since $P_n(z)\mid P_{n+1}(z), \forall n\ge 1$ we will have, for each $k\ge m$, $f_k(\tau)=\b_j$, so $f(\tau)=\b_j$, and thus $f(X)=Y$ and $f^{-1}(Y)=X$.

All the functions $f$ we construct in this manner are transcendental. Indeed, an entire algebraic function which belong to $\Q[[z]]$ is necessarily a polynomial in $\Q[z]$, and all the coefficients $a_n, n\ge 1$ of the Taylor expansion of the functions $f$ we constructed are nonzero, so these functions $f$ cannot be polynomials. Let us also notice that there is an
$\infty$-ary tree of different possibilities for $f$. In fact, if we have chosen $a_1,\ldots, a_{n-1}$, then, in the construction of $f_{n+1}$, there are infinitely many rational possible
choices of $a_n$ (the choices in each step which will depend on $\e_{n,j}$ and $\d_{n,j}$). Thus, we have constructed uncountably
many possible functions $f$ (notice that the set of the algebraic functions $f\in \Q[[z]]$ is a countable set).

%The proof that we can choose $f$ to be transcendental follows because there is an $\infty$-ary tree of different possibilities for $f$ (in each step we have infinitely many possible choices for $\e_{n+1}$, and so for $a_{n+1}$). Thus, we have constructed uncountably many possible functions, and the algebraic entire functions taking $\QQ$ into itself must be polynomials belonging to $\QQ[z]$, which is a countable subset.
\qed

%Acknowledgements

\section*{Acknowledgement}
The authors are grateful to CNPq for the financial support. Part of this work was done during a very enjoyable visit of the first author to IMPA (Rio de Janeiro). He thanks the people of that institution for their hospitality.

% \appendix

% \section{}
% \label{}

%% BIBLIOGRAPHY


\begin{thebibliography}{99}

\bibitem{bookmahler} K. Mahler, \textit{Lectures on Transcendental Numbers}, Lecture Notes in Math., \textbf{546}, Springer-Verlag, Berlin, 1976.

\bibitem{DG} D. Marques, C. G. Moreira, A positive answer for a question proposed by K. Mahler, {\it Math. Ann.} \textbf{368} (2017), 1059--1062.

%\bibitem{crelle} D. Marques, C. G. Moreira, On exceptional sets of transcendental functions with integer coefficients: solution of a Mahler's problem. Submitted.

\bibitem{19} P. St\"{a}ckel, Ueber arithmetische Eingenschaften analytischer Functionen, \textit{Math. Ann.} \textbf{46} (1895), 513--520.

\bibitem{20} P. St\"{a}ckel, Arithmetische eingenschaften analytischer Functionen, \textit{Acta Math.} \textbf{25} (1902), 371--383.

\bibitem{wal1} M. Waldschmidt, Algebraic values of analytic functions, Proceedings of the International Conference on Special Functions and their Applications (Chennai, 2002).  \textit{J. Comput. Appl. Math.} \textbf{160} (2003), 323--333.




\end{thebibliography}
\end{document}